\documentclass{svmult}

\usepackage[latin1]{inputenc}
\usepackage[english]{babel}
\usepackage{amsmath,amssymb}
\usepackage{eufrak}
\usepackage[mathscr]{eucal}
\usepackage{color}
\usepackage[dvips]{graphicx}
\usepackage{fancyhdr}
\usepackage{cite} 
\usepackage{pst-all}
\newcommand{\bbE}{\mathbb{E}}\newcommand{\bbF}{\mathbb{F}}

\newcommand{\bbN}{\mathbb{N}}
\newcommand{\bbP}{\mathbb{P}}\newcommand{\bbR}{\mathbb{R}}

\newcommand{\bbZ}{\mathbb{Z}}

\newcommand{\calB}{\mathcal{B}}

\newcommand{\calP}{\mathcal{P}}
\newcommand{\calU}{\mathcal{U}}

\newcommand{\frakF}{\mathfrak{F}}

\newcommand{\frakP}{\mathfrak{P}}

\newcommand{\frakX}{\mathfrak{X}}

\newcommand{\tonda}[1]{\left( #1 \right)}
\newcommand{\quadra}[1]{\left[ #1 \right]}
\newcommand{\graffa}[1]{\left\{ #1 \right\}}
\newcommand{\norma}[1]{\left\| #1 \right\|}

\newcommand{\modulo}[1]{\left| #1 \right|}
\newcommand{\doubleindex}[2]{
\textrm{\tiny$
\begin{array}{c}
#1 \\
#2
\end{array}
$}}

\newtheorem{Lem}{Lemma}[section]
\newenvironment{lem}[1][]{\begin{Lem}\begin{normalfont}\emph{#1 }}
{\end{normalfont}\finenu\end{Lem}}
\newtheorem{Teo}[Lem]{Theorem}
\newenvironment{teo}[1][]{\begin{Teo}\begin{normalfont}\emph{#1 }}
{\end{normalfont}\finenu\end{Teo}}
\newtheorem{Cor}[Lem]{Corollary}
\newenvironment{cor}[1][]{\begin{Cor}\begin{normalfont}\emph{#1 }}
{\end{normalfont}\finenu\end{Cor}}
\newtheorem{Pro}[Lem]{Proposition}
\newenvironment{pro}[1][]{\begin{Pro}\begin{normalfont}\emph{#1 }}
{\end{normalfont}\finenu\end{Pro}}
\newtheorem{Defi}[Lem]{Definition}
\newenvironment{defi}[1][]{\begin{Defi}\begin{normalfont}\emph{#1 }}
{\end{normalfont}\fine\end{Defi}}
\newtheorem{Oss}[Lem]{Remark}
\newenvironment{oss}[1][]{\begin{Oss}\begin{normalfont}\emph{#1 }}
{\end{normalfont}\fineoss\end{Oss}}
\newtheorem{Es}[Lem]{Example}

\newcounter{HPdisc}[section]
\renewcommand{\theHPdisc}{(Ad-\arabic{HPdisc})}

\newcounter{HPcont}[section]
\renewcommand{\theHPcont}{(A-\arabic{HPcont})}
\newenvironment{assume}{
\begin{itemize}\refstepcounter{HPcont} \item[\!\!\!\!\!\!\!\!\!\!\textbf{\theHPcont}]} {\end{itemize}\par}

\newenvironment{myproof}[1][]{\par\noindent\textbf{Proof{#1}. }}{\finedim\par}

\newcommand{\fine}{}
\newcommand{\finenu}{}
\newcommand{\finedim}{\hfill$\blacksquare$}
\newcommand{\fineoss}{}

\newcounter{numerazione}
\newenvironment{mynumerate}{\setcounter{numerazione}{0}}{}
\newcommand{\mynumber}[1][\refstepcounter{numerazione}(\arabic{numerazione}) ]{\par\noindent #1 }

\newcommand{\myitem}[1][$\bullet$]{\par\noindent #1 $\ \ $}

\newcounter{Bnumero}[Lem]

\newcommand{\Banach}{\frakX} 
\newcommand{\dBanach}{{\Banach^{*}}} 

\newcommand{\unitDBall}{B^{*}_{1}} 

\newcommand{\Mink}{\oplus} 
\newcommand{\Minkowski}{+} 
\newcommand{\HausDist}{\delta_H} 

\newcommand{\Parts}[1]{\frakP(#1)} 
\newcommand{\eParts}[1]{\frakP^{\,0}(#1)}

\newcommand{\Closed}[2][]{\bbF_{#1}(#2)}

\newcommand{\eClosed}[2][]{\bbF^{\,0}_{#1}(#2)}

\newcommand{\closure}[1]{\overline{#1}}
\newcommand{\interior}[1]{\textrm{Int}\ #1}

\newcommand{\nucleation}{B}
\newcommand{\growth}{G}
\newcommand{\cpt}{K}

\newcommand{\salgebra}{\frakF} 
\newcommand{\borel}[1][]{\calB_{#1}}
\newcommand{\prob}[1][]{\bbP_{#1}}

\newcommand{\leb}{\mu_\lambda}
\newcommand{\misura}{\mu}
\newcommand{\expec}{\bbE}

\newcommand{\racs}{{RaCS}}
\newcommand{\grap}{{G-RaP}}

\begin{document}

\title{A set--valued framework for birth--and--growth process}
\author{Giacomo Aletti, Enea G. Bongiorno, Vincenzo Capasso}
\institute{Department of Mathematics, University of Milan,\\ via Saldini 50, 10133 Milan Italy\\
\texttt{giacomo.aletti@mat.unimi.it} \\
\texttt{bongio@mat.unimi.it} \\
\texttt{vincenzo.capasso@mat.unimi.it}}
%
%
\maketitle

\abstract{We propose a set--valued framework for the well--posedness of birth--and--growth process. Our birth--and--growth model is rigorously defined as a suitable combination, involving Minkowski sum and Aumann integral, of two very general set--valued processes representing nucleation and growth respectively. The simplicity of the used geometrical approach leads us to avoid problems arising by an analytical definition of the front growth such as boundary regularities. In this framework, growth is generally anisotropic and, according to a mesoscale point of view, it is not local, i.e.\ for a fixed time instant, growth is the same at each space point.}

\bibliographystyle{plain}

\section*{Introduction}

Nucleation and growth processes arise in several natural and technological applications (cf. \cite{cap03,cap03a} and the references therein) such as, for example, solidification and phase--transition of materials, semiconductor crystal growth, biomineralization, and DNA replication (cf., e.g., \cite{her:jun:bec:ben02}).

A \emph{birth--and--growth process} is a \racs{} family given by $ \Theta_t=\bigcup_{n:T_n\le t} \Theta_{T_n}^t(X_n)$, for $t\in\bbR_{+} $, where $\Theta^t_{T_n}\tonda{X_n}$ is the \racs{} obtained as the evolution up to time $t > T_n$ of the germ born at (random) time $T_n$ in (random) location $X_n$, according to some growth model.
\\%
An analytical approach is often used to model birth--and--growth process, in particular it is assumed that the growth is driven according to a non--negative normal velocity, i.e.\ for every instant $t$, a border point $x\in\partial \Theta_t$ ``grows'' along the outward normal unit (e.g. \cite{bar:sor:sou93, su:bur07, fro:tho87, bur04, bur:cap:piz06, bur:cap:mic07, chi04}). Thus, growth is pointwise isotropic; i.e.\ given a point belonging $\partial\Theta_t$, the growth rate is independently from outward normal direction. Note that, the existence of the outward normal vector imposes a regularity condition on $\partial \Theta_t$ and also on the nucleation process (it cannot be a point process).

This paper is an attempt to offer an original alternative approach based on a purely geometric stochastic point of view, in order to avoid regularity assumptions describing birth--and--growth process. In particular, Minkowski sum (already employed in \cite{mic:pat:vil05} to describe self--similar growth for a single convex germ) and Aumann integral are used here to derive a mathematical model of such process. This model, that emphasizes the geometric growth without regularity assumptions on $\partial\Theta_t$, is rigorously defined as a suitable combination of two very general set--valued processes representing nucleation $\graffa{\nucleation_t}_{t\in[t_0,T]}$ and growth $\graffa{\growth_t}_{t\in[t_0,T]}$ respectively
$$
\begin{array}{rl}
\Theta_t=& \tonda{\Theta_{t_0}\Mink \int_{t_0}^t \growth_s ds}\cup \bigcup_{s\in [t_0, t]} d\nucleation_{s}
\\%
d\Theta_t =& \oplus \growth_t dt \cup d\nucleation_t \qquad\textrm{ or }\qquad
\Theta_{t+dt}=(\Theta_t\Mink \growth_t dt)\cup d\nucleation_t.
\end{array}
$$
Roughly speaking, increment $d\Theta_t$, during an infinitesimal time interval $dt$, is an enlargement due to an infinitesimal Minkowski addend $\growth_{t} dt$ followed by the union with the infinitesimal nucleation $d\nucleation_t$.
\\%
As a consequence of Minkowski sum definition, for every instant $t$, each point $x\in\Theta_t$ (and then each point $x\in\partial \Theta_t$) grows up by $\growth_t dt$ and no regularity border assumptions are required. Then we deal with \emph{not--local} growth; i.e.\ growth is the same Minkowski addend for every $x\in\Theta_t$. Nevertheless, under mesoscale hypothesis we can only consider constant growth region as described, for example, in \cite{bur:cap:piz06}. On the other hand, growth is anisotropic whenever $\growth_t$ is not a ball.


The aim of this paper is to ensure the well--posedness of such a model and, hence, to show that above ``integral'' and ``differential'' notations are meaningful.
\\%
In view of well--posedness, in \cite{ale:bon:cap08b}, the authors show how the model leads to different and significant statistical results.

The article is organized as follows. Section \ref{sec: preliminary results} contains some assumptions about (random) closed sets and their basilar properties.
Model assumptions are collected in Section \ref{sec:assumptions} and integrability properties of growth process are studied in Section \ref{sec:growth_propert}. For the sake of simplicity, we present, in Section \ref{sec: Geometric Random Process}, main results of the paper (that imply well-posedness of the model), whilst correspondent proofs are in Section \ref{sec: proofs of well-position}. At the last, Section \ref{sec:discrete_time_case} proposes a discrete time point of view, also justifying integral and differential notations.

\section{Preliminary results}\label{sec: preliminary results}

Let $\bbN$, $\bbZ$, $\bbR$, $\bbR_{+}$ be the sets of all non--negative integer, integer, real and non--negative real numbers respectively. Let $\Banach$, $\dBanach$, $\unitDBall$ be a Banach space, its dual space and the unit ball of the dual space centered in the origin respectively. We shall consider
$$
\begin{array}{llcll} \eParts{\Banach} &= \textrm{ the family of all
subsets of } \Banach,
& & \Parts{\Banach} &= \eParts{\Banach}\setminus \{\emptyset\}%
\\%
\eClosed{\Banach} &= \textrm{ the family of all closed subsets of
}\Banach, &\quad &\Closed{\Banach} &=
\eClosed{\Banach}\setminus\{\emptyset\}.
\end{array}
$$
The suffixes $c$ and $b$ denote convexity and boundedness properties respectively (e.g. $\eClosed[bc]{\Banach}$ denotes the family of all closed, bounded and convex subsets of $\Banach$).

For all $A,B\in\eParts{\Banach}$ and $\alpha\in\bbR_{+}$, let us
define
$$
\begin{array}{rll}
A\Minkowski B =& \graffa{a+b:a\in A,\ b\in B}=\bigcup_{b\in B}
b\Minkowski A,\quad& \textrm{(Minkowski Sum)}
\\%
\alpha\cdot A=&\alpha A =\graffa{\alpha a: a\in A}, &
\textrm{(Scalar Product)}
\end{array}
$$
By definition, $\forall A\in\eParts{\Banach}$, $\alpha\in\bbR_{+}$, we have $\emptyset\Minkowski A=\emptyset = \alpha \emptyset$.
It is well known that $\Minkowski$ is a commutative and associative operation with a neutral element but $(\Parts{\Banach},\Minkowski)$ is not a group (cf. \cite{rad52}). The following relations are useful in the sequel (see \cite{ser84}): for all $\forall A,B,C\in\Parts{\Banach}$
$$
\begin{array}{c}
(A\cup B)\Minkowski C = (A\Minkowski C)\cup(B\Minkowski C)
\\%
\textrm{if }B\subseteq C, \quad A\Minkowski B\subseteq
A\Minkowski C
\end{array}
$$
In the following, we shall work with closed sets. In general, if $A,B\in\eClosed{\Banach}$ then $A\Minkowski B$ does not belong to $\eClosed{\Banach}$ (e.g., in $\Banach=\bbR$ let $A=\graffa{n+1/n : n> 1}$ and $B=\bbZ$, then $\graffa{1/n=\tonda{n+1/n}+(-n)}\subset A\Minkowski B$ and $1/n\downarrow 0$, but $0\not\in A\Minkowski B$). In view of this fact, we define $A\Mink B =\closure{A\Minkowski B}$ where $\closure{(\cdot)}$ denotes the closure in $\Banach$.

For any $A,B\in\Closed{\Banach}$ the \emph{Hausdorff distance} (or \emph{metric}) is defined by
$$
\HausDist(A,B) = \max\graffa{\sup_{a\in A}\inf_{b\in
B}\norma{a-b}_\Banach,
\sup_{b\in B} \inf_{a\in A}\norma{a-b}_\Banach}.%
$$

For all $(x^*,A)\in \unitDBall\times \Closed{\Banach}$,  the \emph{support function} is defined by $s(x^*,A)=\sup_{ a\in A} x^*(a)$.
It can be proved (cf. \cite{gin:hah:zin83, aub:fra90}) that for each $A,B\in\Closed[bc]{\Banach}$,
\begin{equation}\label{eq: dist Hauss = norm support function}
\HausDist(A,B) = \sup\graffa{ \modulo{ s(x^*,A)-s(x^*,B) } : x^* \in
\unitDBall}.
\end{equation}


Let $(\Omega,\salgebra)$ be a measurable space with $\salgebra$ complete with respect to some $\sigma$-finite measure, let $X:\Omega\to \eParts{\Banach}$ be a set--valued map, and
\begin{align*}
D(X)&=\graffa{\omega\in\Omega: X(\omega)\neq\emptyset} & \textrm{be the \emph{domain} of $X$}
\\%
X^{-1}(A)&=\graffa{\omega\in\Omega : X(\omega)\cap A \neq\emptyset },\quad
A\subset \Banach, & \textrm{be the \emph{inverse image} of $X$}
\end{align*}
Roughly speaking, $X^{-1}(A)$ is the set of all $\omega$ such that $X(\omega)$ hits set $A$.

Different definitions of measurability for set--valued functions are developed over the years by several authors (cf. \cite{him75, cas:val77, aub:fra90, hia:ume77} and reference therein). Here, $X$ is \emph{measurable} if, for each $O$, open subset of $\Banach$, $X^{-1}(O)\in \salgebra$.
\begin{pro}[(See \cite{him75})]%
\label{teo: measurability conditions of RACS}
$X:\Omega\to \eParts{\Banach}$ is a measurable set--valued map if and only if $D(X)\in\salgebra$, and $\omega\mapsto d(x,X(\omega))$ is a measurable
function of $\omega\in D(X)$ for each $x\in\Banach$.
\end{pro}
From now on, $\calU[\Omega,\salgebra,\misura;\Closed{\Banach}]$ ($=\calU[\Omega;\Closed{\Banach}]$ if the measure $\mu$ is clear) denotes the family of $\Closed{\Banach}$--valued measurable maps (analogous notation holds whenever $\Closed{\Banach}$ is replaced by another family of subsets of $\Banach$).

Let $(\Omega,\salgebra,\prob)$ be a complete probability space and let $X\in \calU[\Omega,\salgebra,\prob;\Closed{\Banach}]$, then $X$ is a \racs{}.
\\%
It can be proved (see \cite{li:ogu:kre02}) that, if $X,X_1,X_2$ are \racs{} and if $\xi$ is a measurable real--valued function, then $X_1\oplus X_2$, $X_1\ominus X_2$, $\xi X$ and $(\interior{X})^C$ are \racs{}. Moreover, if $\graffa{X_n}_{n\in\bbN}$ is a sequence of \racs{} then $X=\closure{\bigcup_{n\in\bbN} X_n}$ is so.

Let $(\Omega,\salgebra,\misura)$ be a finite measure space (although most of the results are valid for $\sigma$-finite measures space). The \emph{Aumann integral} of $X\in \calU[\Omega,\salgebra,\misura;\Closed{\Banach}]$ is defined by
$$
\int_\Omega Xd\misura=\graffa{\int_\Omega x d\misura:x\in S_X},
$$
where $S_X=\graffa{x\in L^1[\Omega;\Banach]: x\in X\ \mu-\textrm{a.e.}}$ and $\int_\Omega x d\misura$ is the usual Bochner integral in $L^1[\Omega;\Banach]$. Moreover, $\int_A X d\misura=\graffa{\int_A x d\misura:x\in S_X}$ for $A\in\salgebra$. If $\misura$ is a probability measure, we denote the Aumann integral by $\expec{X} = \int_\Omega Xd\misura$.

Let $X\in\calU[\Omega,\salgebra,\misura;\Closed{\Banach}]$, it is \emph{integrably bounded}, and we shall write $X\in L^1[\Omega,\salgebra,\misura;\Closed{\Banach}]=L^1[\Omega;\Closed{\Banach}]$, if $\norma{X}_h\in L^1[\Omega,\salgebra,\misura;\bbR]$.

\section{Model assumptions}
\label{sec:assumptions}

Let us consider
\begin{equation}\label{eq: infinitesimal continuous like discrete}
\begin{array}{rl}
\Theta_t=& \tonda{\Theta_{t_0}\Mink \int_{t_0}^t \growth_s ds}\cup \bigcup_{s\in [t_0, t]} d\nucleation_{s}
\\%
d\Theta_t =& \oplus \growth_t dt \cup d\nucleation_t \qquad\textrm{ or }\qquad
\Theta_{t+dt}=(\Theta_t\Mink \growth_t dt)\cup d\nucleation_t.
\end{array}
\end{equation}
In fact, above equation is not a definition since, for example, problems arise handling non--countable union of (random) closed sets. The well--posedness of \eqref{eq: infinitesimal continuous like discrete} and hence the existence of such a process are the main purpose of this paper.
\\%
From now on, let us consider the following assumptions.
\addtocounter{HPcont}{-1}
\begin{assume}\label{hp(cont): basic hypothesis}
\begin{itemize}
\item[-]
$(\Banach,\norma{\cdot}_\Banach )$ is a reflexive Banach space with separable dual space $(\dBanach,\norma{\cdot}_{\dBanach})$, (then, $\Banach$ is separable too, see \cite[Lemma II.3.16 p. 65]{dun:sch58}).

\item[-]
$[t_0,T]\subset \bbR $ is the \emph{time observation interval} (or
\emph{time interval}),

\item[-]
$\tonda{ \Omega,\salgebra,\graffa{\salgebra_t}_{t\in [t_0,T]}, \prob }$ is a filtered probability space, where the filtration $\graffa{\salgebra_t}_{t\in [t_0,T]}$ is assumed to have the usual properties.
\end{itemize}
\end{assume}

\noindent(\emph{Nucleation Process}). $\nucleation=\graffa{\nucleation(\omega,t)=\nucleation_t:\omega\in\Omega,\ t\in[t_0,T]}$ is a process with non--empty closed values, i.e.\
$\nucleation: \Omega \times [t_0,T] \to \Closed{\Banach}$
such that
\begin{assume}\label{hp(cont): Ht RACS}
$\nucleation(\cdot,t)\in\calU[\Omega,\salgebra_t,\prob; \Closed{\Banach}]$, for every $t\in [t_0,T]$, i.e. $\nucleation_t$ is an \emph{adapted} (to $\graffa{\salgebra_t}_{t\in [t_0,T]}$) process.
\end{assume}
\begin{assume}\label{hp(cont): Ht increasing}
$\nucleation_t$ is increasing: for every $t,s\in [t_0,T]$ with $s< t$, $\nucleation_s
\subseteq \nucleation_{t}$.
\end{assume}

\noindent(\emph{Growth Process}). $\growth=\graffa{\growth_t=\growth(\omega,t):\omega\in\Omega,\ t\in[t_0,T]
}$ is a process with non--empty closed values, i.e.\
$\growth: \Omega \times [t_0,T] \to \Closed{\Banach}$ such that
\begin{assume}\label{hp(cont): 0 belongs to G}
for every $\omega\in\Omega$ and $t\in[t_0,T]$, $0\in\growth(\omega,t)$.
\end{assume}
\begin{assume}\label{hp(cont): G convex}
for every $\omega\in\Omega$ and $t\in[t_0,T]$, $\growth(\omega,t)$ is
convex, i.e.\ $\growth: \Omega \times [t_0,T] \to \Closed[c]{\Banach}$.
\end{assume}
\begin{assume}\label{hp(cont): G bounded by K}
there exists $\cpt\in\Closed[b]{\Banach}$ such that $\growth(\omega,
t)\subseteq \cpt$ for every $t\in[t_0,T]$ and $\omega\in\Omega$.
\end{assume}
\noindent As a consequence, $\growth(\omega,t) \in \Closed[b]{\Banach}$ and $\norma{\growth(\omega,t)}_h\le\norma{\cpt}_h$, $\forall (\omega,t)\in \Omega \times [t_0,T]$.
%
%

In order to establish the well--posedness of integral $\int_{t_0}^t \growth_s ds$ in \eqref{eq: infinitesimal continuous like discrete}, let us consider a suitable hypothesis of measurability for $\growth$ (analogously to what is).
\\%
A $\Closed{\Banach}$--valued process $\growth=\graffa{\growth_t}_{t\in [t_0,T]}$  has \emph{left continuous trajectories} on $[t_0,T]$ if, for every $\overline{t}\in [t_0,T]$ with $t<\overline{t}$,
$$
\lim_{t\to \overline{t}} \HausDist\tonda{\growth(\omega,t) ,
\growth(\omega,\overline{t})} =0,\qquad \textrm{a.s.}
$$
\\%
The $\sigma$-algebra on $\Omega\times[t_0,T]$ generated by the processes
$\graffa{\growth_t}_{t\in[t_0,T]}$ with left continuous trajectories on $[t_0,T]$, is called the \emph{previsible} (or \emph{predictable}) $\sigma$-algebra and it is denoted by $\calP$.
\begin{pro}
The previsible $\sigma$-algebra is also generated by the collection of
random sets $A\times t_0$ where $A\in\salgebra_{t_0}$ and $A\times (s,t]$
where $A\in\salgebra_s$ and $(s,t]\subset[t_0,T]$.
\end{pro}
\begin{myproof}
Let the $\sigma$-algebra generated by the above collection of sets be
denoted by $\calP'$. We shall show $\calP=\calP'$. Let $\growth$ be a left
continuous process and let $\alpha=(T-t_0)$, consider for $n\in\bbN$
$$
\growth_n(\omega,t)=\left\{
\begin{array}{ll}
\growth(\omega,t_0), & t=t_0\\
\\
\growth\tonda{\omega,t_0+\frac{k\alpha}{2^n}},\  &
\begin{array}{c}
\tonda{t_0+\frac{k\alpha}{2^n}} < t \le
\tonda{t_0+\frac{(k+1)\alpha}{2^n}}
\\
k\in\graffa{0,\ldots, (2^n-1)}
\end{array}
\end{array}\right.
$$
It is clear that $\growth_n$ is $\calP'$-measurable, since $\growth$ is
adapted. As $\growth$ is left continuous, the above sequence of
left-continuous processes converges pointwise (with respect to $\HausDist$) to $\growth$ when $n$ tends to infinity, so $\growth$ is
$\calP'$-measurable, thus $\calP\subseteq\calP'$.
\\%
Conversely consider $A\times(s,t]\in\calP'$ with $(s,t]\subset[t_0,T]$ and
$A\in\salgebra_s$. Let $b\in\Banach\setminus\{0\}$ and $\growth$ be the
process
$$
\growth(\omega,v)=\left\{
\begin{array}{ll}
b, & v\in (s,t],\ \omega\in A\\
0, & \textrm{otherwise}
\end{array}\right.
$$
this function is adapted and left continuous, hence
$\calP'\subseteq\calP$.
\end{myproof}

Then let us consider the following assumption.
\begin{assume}\label{hp(cont): G is P-measurable}
$\growth$ is $\calP$-measurable.
\end{assume}

\section{Growth process properties}%
\label{sec:growth_propert}

Theorem \ref{teo: any lebesg-integral of G is RACS} is the main result in this section. It shows that $\omega\mapsto \int_a^b \growth(\omega,\tau) d\tau$ is a \racs{} with non--empty bounded convex values. This is the first step in order to obtain well--posedness of \eqref{eq: infinitesimal continuous like discrete}.

%
\begin{pro}\label{pro: G(w,t) integrably bounded}
Suppose \ref{hp(cont): 0 belongs to G}, \ldots, \ref{hp(cont): G is P-measurable}, and let $\leb$ be the Lebesgue measure on $[t_0,T]$, then
\begin{itemize}
\item
$\growth(\omega,\cdot)\in
\calU\quadra{[t_0,T],\borel[{[t_0,T]}],\leb;\Closed[bc]{\Banach}}$ for
every $\omega\in\Omega$.

\item
$\growth(\cdot,t)\in \calU[\Omega,\widetilde{\salgebra}_{t^-}, \prob;
\Closed[bc]{\Banach}]$ for each $t\in[t_0,T]$, where
$\widetilde{\salgebra}_{t^-}$ is the so called \emph{history
$\sigma$-algebra} i.e.\
$\widetilde{\salgebra}_{t^-}=\sigma\tonda{\salgebra_s : 0\le s <
t}\subseteq\salgebra$.

\item
$\growth\in\
L^1[[t_0,T], \borel[{[t_0,T]}], \leb; \Closed[bc]{\Banach}]
\cap L^1[\Omega,\salgebra,\prob; \Closed[bc]{\Banach}]$
\end{itemize}
\end{pro}
\begin{myproof}
Assumptions \ref{hp(cont): 0 belongs to G} and \ref{hp(cont): G convex} imply that $\growth$ is non--empty and convex. Measurability and  integrability properties are consequence of \ref{hp(cont): G is P-measurable} and \ref{hp(cont): G bounded by K} respectively.
\end{myproof}
\begin{teo}\label{teo: any lebesg-integral of G is RACS}
Suppose \ref{hp(cont): 0 belongs to G}, \ldots, \ref{hp(cont): G is P-measurable}. For every $a, b\in [t_0,T]$, the integral $\int_a^b \growth(\omega,\tau) d\tau$ is non--empty and the set--valued map
$$
\begin{array}{rccl}
\growth_{a,b}: & \Omega &\to& \Parts{\Banach}
\\
&\omega &\mapsto& \int_a^b \growth(\omega,\tau) d\tau
\end{array}
$$
is measurable. Moreover, $\growth_{a,b}$ is a non--empty, bounded convex \racs.
\end{teo}

In order to prove Theorem \ref{teo: any lebesg-integral of G is RACS}, consider following properties for real processes.
\\%
A real--valued process $X=\graffa{X_t}_{t\in [t_0,T]}$ is \emph{predictable} with respect to filtration $\graffa{\salgebra_t}_{t\in\bbR_{+}}$, if it is measurable with respect to the \emph{predictable $\sigma$-algebra} $\calP_{\bbR}$, i.e.\ the $\sigma$-algebra generated by the collection of random sets $A\times\graffa{0}$ where $A\in\salgebra_0$ and $A\times (s,t]$ where $A\in\salgebra_s$.
\begin{pro}[{(See \cite[Propositions 2.30, 2.32 and 2.41]{cap:bak05})}]%
\label{pro: real-valued stochastic results}%
Let $X=\graffa{X_t}_{t\in [t_0,T]}$ be a predictable real--valued process, then $X$ is \makebox{$(\salgebra\otimes\borel[{[t_0,T]}],\borel[\bbR])$}-measurable.
Further, for every $\omega \in \Omega$, the trajectory $X(\omega,\cdot):
[t_0,T]\to \bbR$ is $(\borel[{[t_0,T]}],\borel[\bbR])$-measurable .
\end{pro}
\begin{lem}\label{lem: support function measurable}
Let $x^*$ be an element of the unit ball in the dual space $\unitDBall$,
then $\growth \mapsto s(x^*,\growth)$ is a measurable map.
\end{lem}
\begin{myproof}
By definition $s(x^*,\growth)=\sup\graffa{ x^*(g):g\in\growth }$.
Since $\Banach$ is separable \ref{hp(cont): basic hypothesis}, there exists $\graffa{g_n}_{n\in\bbN}\subset \growth$ such that $\growth=\overline{ \graffa{g_n} }$. Then, for every
$x^*\in \unitDBall$ we have
$$
s(x^*,\growth)=\sup_{g\in\growth} x^*(g)=\sup_{n\in \bbN} x^*(g_n).
$$
Since $x^*$ is a continuous map then, $s(x^*,\cdot)$ is measurable.
\end{myproof}
\begin{myproof}[ of Theorem \ref{teo: any lebesg-integral of G is RACS}]
At first, we prove that $\growth_{a,b}$ is a measurable map. From Proposition \ref{pro: G(w,t) integrably bounded}, integral $ \growth_{a,b}=\int_a^b \growth(\omega,\tau) d\tau $ is well defined for all $\omega\in\Omega$. Assumption \ref{hp(cont): 0 belongs to G} implies $0\in \growth_{a,b}(\omega)\neq\emptyset$ for every $\omega\in\Omega$. Hence, the domain of $\growth_{a,b}$ is the whole $\Omega$ for all $a,b\in[t_0,T]$
$$
D\tonda{\growth_{a,b}}=\graffa{\omega\in\Omega: \growth_{a,b}
\neq\emptyset}=\Omega\in\salgebra.
$$
Thus, by Proposition \ref{teo: measurability conditions of RACS} and for a
fixed couple $a,b\in[t_0,T]$, $\growth_{a,b}$ is (weakly) measurable if and
only if, for every $x\in\Banach$, the map
\begin{equation}\label{eq: 1inProp.any lebesg-integral of G is RACS}
\omega \mapsto d\tonda{x,\int_a^b\growth(\omega,\tau)d\tau}=
\HausDist\tonda{x,\int_a^b\growth(\omega,\tau)d\tau}
\end{equation}
is measurable. Equation \eqref{eq: dist Hauss = norm support function} guarantees that \eqref{eq: 1inProp.any lebesg-integral of G is RACS} is measurable if and only if, for every $x\in\Banach$, the map
$$
\omega\mapsto \sup_{x^*\in \unitDBall} \modulo{
s(x^*,x)-s\tonda{x^*,\int_a^b\growth(\omega,\tau)d\tau} }
$$
is measurable. The above expression can be computed on a countable family
dense in $\unitDBall$ (note that such family exists since $\dBanach$ is
assumed separable \makebox{\ref{hp(cont): basic hypothesis}}):
$$
\omega\mapsto \sup_{n\in \bbN} \modulo{
s(x^*_i,x)-s\tonda{x^*_i,\int_a^b\growth(\omega,\tau)d\tau} }.
$$
It can be proved (\cite[Theorem 2.1.12 p. 46]{li:ogu:kre02}) that
$$
s\tonda{x^*,\int_a^b\growth(\omega,\tau)d\tau}= \int_a^b
s\tonda{x^*,\growth(\omega,\tau)}d\tau, \qquad\forall x^*\in\unitDBall
$$
and therefore, since $s(x^*_i,x)$ is a constant, $\growth_{a,b}$ is
measurable if, for every $x^*\in \graffa{x_i^*}_{i\in\bbN}$, the following
map
\begin{equation}\label{eq: 2inProp.any lebesg-integral of G is RACS}
\begin{array}{ccl}
(\Omega,\salgebra) & \to & (\bbR,\borel[\bbR])
\\%
\omega&\mapsto &\int_a^b s\tonda{x^*,\growth(\omega,\tau)}d\tau
\end{array}
\end{equation}
is measurable. Note that $s(x^*,\growth(\cdot,\cdot))$, as a map from
$\Omega\times[t_0,T]$ to $\bbR$, is predictable since it is the
composition of a predictable map \ref{hp(cont): G is P-measurable} with a
measurable one (see Lemma \ref{lem: support function measurable}):
$$
\begin{array}{rccccc}
s\tonda{x^*,\growth(\cdot,\cdot)}:& (\Omega\times[t_0,T],\calP) &\to&
(\Closed{\Banach},\sigma_f) &\to& (\bbR,\borel[\bbR])
\\
& (\omega,t) &\mapsto& \growth(\omega,t) &\mapsto&
s\tonda{x^*,\growth(\omega,t)}
\end{array}
$$
thus, by Proposition \ref{pro: real-valued stochastic results}, it is a
$\calP$-measurable map and hence \eqref{eq: 2inProp.any lebesg-integral of
G is RACS} is a measurable map.

In view of the first part, it remains to prove that $\growth_{a,b}$
is a bounded convex set for a.e. $\omega\in\Omega$. Since $\Banach$ is reflexive \ref{hp(cont): basic hypothesis}, by Proposition \ref{pro: G(w,t) integrably bounded} we have that $\growth_{a,b}$ is closed (\cite[Theorem 2.2.3]{li:ogu:kre02}).
Further, $\growth_{a,b}$ is also convex (see \cite[Theorem 2.1.5 and Corollary 2.1.6]{li:ogu:kre02}).
\\%
To conclude the proof, it is sufficient to show that $\growth_{a,b}$ is included in a bounded set:
\begin{eqnarray*}
\int_a^b \growth(\omega,\tau)d\tau &=& \graffa{\int_a^b g(\omega,\tau)
d\tau :
g(\omega,\cdot)\in \growth(\omega,\cdot)\subseteq \cpt}\\
&\subseteq& \graffa{\int_a^b k d\tau : k\in \cpt}= \graffa{ (b-a)k  : k\in \cpt} = (b-a)\cpt.
\end{eqnarray*}
\end{myproof}

\section{Geometric Random Process}%
\label{sec: Geometric Random Process}%

For the sake of simplicity, let us present the main results which proofs will be given in Section \ref{sec: proofs of well-position}.
\\%
Let us assume conditions from \ref{hp(cont): basic hypothesis} to
\ref{hp(cont): G is P-measurable}. For every $t\in[t_0,T]\subset\bbR$, $n\in\bbN$ and $\Pi=\tonda{t_i}_{i=0}^n$ partition of $[t_0,t]$, let us define
\begin{align}
s_{\Pi}(t) = &\tonda{\nucleation_{t_0} \Mink \int_{t_0}^t \growth(\tau) d\tau
}\cup \bigcup_{i=1}^{n} \tonda{\Delta \nucleation_{t_i}\Mink \int_{t_i}^t
\growth(\tau) d\tau}\label{eq: lower union}
\\
S_{\Pi}(t) = &\tonda{\nucleation_{t_0} \Mink \int_{t_0}^t \growth(\tau) d\tau
}\cup \bigcup_{i=1}^{n} \tonda{\Delta \nucleation_{t_i}\Mink \int_{t_{i-1}}^t
\growth(\tau) d\tau}
\label{eq: upper union}
\end{align}
where $\Delta\nucleation_{t_i} = \nucleation_{t_i}\setminus \nucleation_{t_{i-1}}^o$
($\nucleation_{t_{i-1}}^o$ denotes the interior set of $\nucleation_{t_{i-1}}$) and where the integral is in the Aumann sense with respect to the Lebesgue measure $d\tau = d\leb$. We write $s_{\Pi}$ and $S_\Pi$ instead of $s_\Pi(t)$ and $S_\Pi (t)$ when the dependence on $t$ is clear.

Proposition \ref{pro: sn and Sn are RACS} guarantees that both $s_{\Pi}$ and $S_\Pi$ are well defined \racs{}, further, Proposition \ref{pro: sn<Sn} shows $s_\Pi\subseteq S_\Pi$ as a consequence of different time intervals integration: if the time interval integration of $\growth$ increases then the integral of $\growth$ does not decrease with respect to set-inclusion (Lemma \ref{lem: int_i G < int_I G}). Proposition \ref{pro: sn<sm and Sm<Sn} means that $\graffa{s_\Pi}$ ($\graffa{S_\Pi}$) increases (decreases) whenever a refinement of $\Pi$ is considered. At the same time, Proposition \ref{pro: dist(sn,Sn)--> 0} implies that $s_\Pi$ and $S_\Pi$ become closer each other (in the Hausdorff distance sense) when partition $\Pi$ becomes finer. The ``limit'' is independent on the choice of the refinement as consequence of Proposition \ref{pro: independence of partition}.
\\%
Corollary \ref{cor: definition of theta} means that, given any
$\graffa{\Pi_j}_{j\in\bbN}$ refinement sequence of $[t_0,t]$, the random
closed sets $s_{\Pi_j}$ and $S_{\Pi_j}$ play the same role that lower sums
and upper sums have in classical analysis when we define the Riemann
integral. In fact, if $\Theta_t$ denotes their limit value (see \eqref{eq:
limits sn Sn}), $s_{\Pi_j}$ and $S_{\Pi_j}$ are a lower and an upper
approximation of $\Theta_t$ respectively.
Note that, as a consequence of monotonicity of $s_{\Pi_j}$ and $S_{\Pi_j}$, we avoid problems that may arise considering uncountable unions in integral expression in \eqref{eq: infinitesimal continuous like discrete}.

\begin{pro}\label{pro: sn and Sn are RACS}
Let $\Pi$ be a partition of $[t_0,t]$. Both $s_\Pi$ and $S_\Pi$, defined
in \eqref{eq: lower union} and \eqref{eq: upper union}, are \racs{}.
\end{pro}
\begin{lem}\label{lem: int_i G < int_I G}
Let $X\in L^1[I,\salgebra,\leb;\Closed{\Banach}]$, where $I$ is a bounded
interval of $\bbR$, such that $0\in X$ $\leb$-almost everywhere on $I$ and
let $I_1,I_2$ be two other intervals of $\bbR$ with $I_1\subset I_2\subset
I$. Then
$$
\int_{I_1} X(\tau) d\tau \subseteq \int_{I_2} X(\tau) d\tau.
$$
\end{lem}
\begin{pro}\label{pro: sn<Sn}
Let $\Pi$ be a partition of $[t_0,t]$. Then $s_\Pi\subseteq S_\Pi$ almost
surely.
\end{pro}
\begin{pro}\label{pro: sn<sm and Sm<Sn}
Let $\Pi$ and $\Pi'$ be two partitions of $[t_0,t]$ such that $\Pi'$ is a
refinement of $\Pi$. Then, almost surely, $s_\Pi\subseteq s_{\Pi'}$ and
$S_{\Pi'}\subseteq S_\Pi$.
\end{pro}
\begin{pro}\label{pro: dist(sn,Sn)--> 0}
Let $\graffa{\Pi_j}_{j\in\bbN}$ be a refinement sequence of $[t_0,t]$ (i.e.\ $\modulo{\Pi_j}\to 0$ if $j\to\infty$). Then, almost surely,
$\lim_{j\to\infty}\HausDist\tonda{s_{\Pi_j},S_{\Pi_j}} = 0$.
\end{pro}
\begin{pro}\label{pro: independence of partition}
Let $\graffa{\Pi_j}_{j\in\bbN}$ and $\graffa{\Pi_l'}_{l\in\bbN}$ be two
distinct refinement sequences of $[t_0,t]$, then, almost surely,
$$
\lim_\doubleindex{j\rightarrow \infty}{l\rightarrow \infty}
\HausDist\tonda{s_{\Pi_j},s_{\Pi'_l}} = 0 \qquad \textrm{and}\qquad
\lim_\doubleindex{j\rightarrow \infty}{l\rightarrow \infty}
\HausDist\tonda{S_{\Pi_j},S_{\Pi'_l}} = 0.
$$
\end{pro}
\begin{cor}\label{cor: definition of theta}
For every $\graffa{\Pi_j}_{j\in\bbN}$ refinement sequence of $[t_0,t]$,
the following limits exist
\begin{equation}\label{eq: limits sn Sn}
\overline{\tonda{\bigcup_{j\in\bbN} s_{\Pi_j}}},\
\overline{\tonda{\lim_{j\rightarrow\infty} s_{\Pi_j}}},\
\lim_{j\rightarrow \infty} S_{\Pi_j},\ \bigcap_{j\in\bbN} S_{\Pi_j},
\end{equation}
and they are equals almost surely. The convergences is taken with respect
to the Hausdorff distance.
\end{cor}
We are now ready to define the continuous time stochastic process.
\begin{defi}\label{def: continuous set process}
Assume \ref{hp(cont): basic hypothesis}, \ldots, \ref{hp(cont): G is
P-measurable}. For every $t\in [t_0,T]$, let $\graffa{\Pi_j}_{j\in\bbN}$
be a refinement sequence of the time interval $[t_0,t]$ and let $\Theta_t$
be the \racs{} defined by
$$
\overline{\tonda{\bigcup_{j\in\bbN} s_{\Pi_j}(t)}}=
\overline{\tonda{\lim_{j\rightarrow\infty} s_{\Pi_j}(t)}} = \Theta_t=
\lim_{j\rightarrow \infty} S_{\Pi_j}(t) = \bigcap_{j\in\bbN} S_{\Pi_j}(t),
$$
then, the family $\Theta=\graffa{\Theta_t : t\in [t_0,T]}$ is called
\emph{geometric random process \grap{}} (on $[t_0,T]$).

%
\end{defi}
\begin{teo}\label{teo: theta increasing}
Let $\Theta$ be a \grap{} on $[t_0,T]$, then $\Theta$ is a non-decreasing process with respect to the set inclusion, i.e.\
$$
\prob\tonda{\Theta_s\subseteq\Theta_t,\ \forall t_0\le s<t\le T}=1.
$$
Moreover, $\Theta$ is adapted with respect to filtration $\graffa{\salgebra_t}_{t\in[t_0,T]}$.
\end{teo}
\begin{oss}
We want to point out that, assumptions we considered on $\graffa{\nucleation_t}$ and $\graffa{\growth_t}$ are so general, that a wide family of classical random sets and evolution processes can be described (for example, Boolean model is a birth--and--growth process with ``null growth'').
\end{oss}

\subsection{Proofs of Propositions in Section \ref{sec: Geometric Random Process}}
\label{sec: proofs of well-position}%

\begin{myproof}[ of Proposition \ref{pro: sn and Sn are RACS}]
For every $i\in\graffa{0,\ldots, n}$, $\int_{t_{i-1}}^t \growth(\tau) d\tau$ is a \racs{} (Theorem \ref{teo: any lebesg-integral of G is RACS}). Thus, measurability Assumption \ref{hp(cont): Ht RACS} on $\nucleation$ guarantees that, for every $t_i \in\Pi$,
$ \nucleation_{t_i}$, $\Delta \nucleation_{t_i}$, $\tonda{\Delta \nucleation_{t_i}\Mink \int_{t_{i}}^t \growth(\tau) d\tau}$, and hence $s_\Pi$ and $S_\Pi$ are \racs{}.
\end{myproof}
\begin{myproof}[ of Lemma \ref{lem: int_i G < int_I G}]
Let $y\in\tonda{\int_{I_1} X(\tau) d\tau}$, then there exists $x\in S_X$, for which $y=\tonda{\int_{I_1} x(\tau) d\tau}$.
Let us define on $I_2(\supset I_1)$
$$
x'(\tau)=\left\{
\begin{array}{ll}
x(\tau), & \tau\in I_1
\\%
0, & \tau\in I_2\setminus I_1
\end{array}
\right.
$$
then $x'\in S_X$ and $y=\tonda{\int_{I_2} x'(\tau) d\tau}\in
\tonda{\int_{I_2} X(\tau) d\tau}$.
\end{myproof}
\begin{myproof}[ of Proposition \ref{pro: sn<Sn}]
Thesis is a consequence of Lemma \ref{lem: int_i G < int_I G}
and Minkowski addition properties, in fact
$\tonda{\int_{t_{i-1}}^t \growth(\tau) d\tau } \subseteq
\tonda{\int_{t_i}^t \growth(\tau) d\tau }$ implies $s_\Pi\subseteq S_\Pi$.
\end{myproof}
\begin{myproof}[ of Proposition \ref{pro: sn<sm and Sm<Sn}]
Let $\Pi'$ be a refinement of partition $\Pi$ of $[t_0,t]$, i.e.\
$\Pi\subset \Pi'$. We prove that $s_\Pi\subseteq s_{\Pi'}$ ($S_{\Pi'}\subseteq S_\Pi$ is analogous). It is sufficient to show the thesis only for $\Pi' = \Pi\cup\graffa{\overline{t}}$ where $\Pi=\graffa{t_0,\ldots, t_n}$ with $t_0< \ldots< t_n=t$ and $\overline{t}\in(t_0,t)$. Let $i\in\graffa{0,\ldots,(n-1)}$ be such that $t_i\le\overline{t}\le t_{i+1}$
then
\begin{eqnarray*}
s_\Pi &=& \tonda{\nucleation_{t_0} \Mink \int_{t_0}^t \growth(\tau) d\tau
} \cup \bigcup_\doubleindex{j=1}{j\neq
i+1}^n\tonda{\Delta\nucleation_{t_j}\Mink\int_{t_j}^t\growth(\tau)d\tau} \cup
\\
&&\quadra{\tonda{\nucleation_{t_{i+1}}\setminus\nucleation_{t_i}^o
}\Mink\int_{t_{i+1}}^t\growth(\tau)d\tau}
\end{eqnarray*}
and
\begin{eqnarray*}
s_{\Pi'} &=& \tonda{\nucleation_{t_0} \Mink \int_{t_0}^t \growth(\tau) d\tau
} \cup \bigcup_\doubleindex{j=1}{j\neq
i+1}^n\tonda{\Delta\nucleation_{t_j}\Mink\int_{t_j}^t\growth(\tau)d\tau} \cup
\\
&&\quadra{\tonda{\nucleation_{\overline{t}}\setminus\nucleation_{t_i}^o
}\Mink\int_{\overline{t}}^t\growth(\tau)d\tau} \cup
\quadra{\tonda{\nucleation_{t_{i+1}}\setminus\nucleation_{\overline{t}}^o
}\Mink\int_{t_{i+1}}^t\growth(\tau)d\tau}
\end{eqnarray*}
Definitely, in order to prove that $s_\Pi\subseteq s_{\Pi'}$ we have to
prove that
\begin{eqnarray*}
\graffa{\quadra{\tonda{\nucleation_{\overline{t}}\setminus\nucleation_{t_i}^o
}\Mink\int_{\overline{t}}^t\growth(\tau)d\tau} \cup
\quadra{\tonda{\nucleation_{t_{i+1}}\setminus\nucleation_{\overline{t}}^o
}\Mink\int_{t_{i+1}}^t\growth(\tau)d\tau}} \\ \supseteq
\quadra{\tonda{\nucleation_{t_{i+1}}\setminus\nucleation_{t_i}^o
}\Mink\int_{t_{i+1}}^t\growth(\tau)d\tau}.
\end{eqnarray*}
This inclusion is a consequence of
$\tonda{\int_{\overline{t}}^t\growth(\tau)d\tau}\supseteq
\tonda{\int_{t_{i+1}}^t\growth(\tau)d\tau}$ (Lemma \ref{lem: int_i G <
int_I G}) and of the Minkowski distribution property.
\end{myproof}
\begin{myproof}[ of Proposition \ref{pro: dist(sn,Sn)--> 0}]
Let $\Pi_j=\tonda{t_i}_{i=0}^n$ be the $j$-partition of the refinement
sequence $\graffa{\Pi_j}_{j\in\bbN}$, then
$$
\HausDist\tonda{s_{\Pi_j},S_{\Pi_j}} = \max\graffa{ \sup_{x\in s_{\Pi_j}}
d(x,S_{\Pi_j}) , \sup_{y\in S_{\Pi_j}} d(y,s_{\Pi_j})}
$$
where $d(x,S_{\Pi_j})=\inf_{y\in S_{\Pi_j}}\norma{x-y}_\Banach$. By
Proposition \ref{pro: sn<Sn}, $s_{\Pi_j}\subseteq S_{\Pi_j}$ then
$$
\sup_{x\in s_{\Pi_j}} d(x,S_{\Pi_j}) = 0
$$
and hence we have to prove that, whenever $j\to \infty$ (i.e.\ $\modulo{\Pi_j}\to 0$),
$$
\HausDist\tonda{s_{\Pi_j},S_{\Pi_j}} = \sup_{y\in S_{\Pi_j}}
d(y,s_{\Pi_j}) = \sup_{y\in S_{\Pi_j}} \inf_{x\in
s_{\Pi_j}}\norma{x-y}_\Banach \longrightarrow 0.
$$

For every $\omega\in\Omega$, let $y$ be any element of
$S_{\Pi_j}(\omega)$, then we distinguish two cases:
\begin{mynumerate}
\mynumber if $y\in \tonda{\nucleation_{t_0}(\omega) \Mink \int_{t_0}^t
\growth(\omega,\tau) d\tau }$, then it is also an element of
$s_{\Pi_j}(\omega)$, and hence $d\tonda{s_{\Pi_j}(\omega),y}=0$.

\mynumber if $y\not\in \tonda{\nucleation_{t_0}(\omega) \Mink \int_{t_0}^t
\growth(\omega,\tau) d\tau }$, then there exist $j\in\graffa{1,\ldots,n}$
such that
$$
y\in \tonda{\Delta\nucleation_{t_j}(\omega)\Mink \int_{t_{j-1}}^t
\growth(\omega,\tau)d\tau }.
$$
By definition of $\Mink$, for every $\omega\in\Omega$, there exist
$$
\graffa{y_m}_{m\in\bbN}\subseteq \tonda{\Delta\nucleation_{t_j}(\omega)
\Minkowski \int_{t_{j-1}}^t \growth(\omega,\tau)d\tau },
$$
such that $\lim_{m\to\infty}y_m=y$. Then, for every $\omega\in\Omega$,
there exist $h_m\in \Delta\nucleation_{t_j}(\omega)$ and $g_m\in
\tonda{\int_{t_{j-1}}^t \growth(\omega,\tau)d\tau}$ such that
$y_m=(h_m+g_m)$ and hence
$$
y=\lim_{m\to\infty} (h_m + g_m)=\lim_{m\to\infty} y_m
$$
where the convergence is in the Banach norm, then let
$\overline{m}\in\bbN$ be such that $\norma{y-y_m}_{\Banach} <
\modulo{\Pi_j}$, for every $m>\overline{m}$.

Note that, for every $\omega\in\Omega$ and $m\in\bbN$, by Aumann integral
definition, there exists a selection
$\widehat{g_m}(\cdot)$ of $ \growth(\omega,\cdot)$ (i.e.\
$\widehat{g_m}(t)\in\growth(\omega,t)$ $\leb$-a.e.) such that
$$
g_m = \int_{t_{j-1}}^t \widehat{g_m}(\tau)d\tau\qquad \textrm{ and }\qquad
y_m = {h_m + \int_{t_{j-1}}^t \widehat{g_m}(\tau)d\tau}.
$$
For every $\omega\in\Omega$, let us consider
$$
x_m= h_m + \int_{t_{j}}^{t} \widehat{g_m}(\tau) d\tau
$$
then $x_m\in s_{\Pi_j}(\omega)$ for all $m\in\bbN$. Moreover, the
following chain of inequalities hold, for all $m > \overline{m}$ and
$\omega\in\Omega$,
\begin{eqnarray*}
\inf_{x'\in s_{\Pi_j}}\norma{x'-y}_{\Banach} &\le& \norma{x_m-y}_{\Banach}
\le \norma{x_m -y_m}_{\Banach} + \norma{y_m - y}_{\Banach}
\\%
&\le& \norma{\int_{t_{j-1}}^{t_j} \widehat{g_m}(\tau) d\tau}_{\Banach}  +
\modulo{\Pi_j} \le \int_{t_{j-1}}^{t_j}
\norma{\widehat{g_m}(\tau)}_{\Banach} d\tau  + \modulo{\Pi_j}
\\%
&\le& \int_{t_{j-1}}^{t_j} \norma{\growth(\tau)}_h d\tau  + \modulo{\Pi_j}
\le \modulo{t_j - t_{j-1}}\norma{\cpt}_h + \modulo{\Pi_j}
\\%
&\le&  \modulo{\Pi_j} \tonda{\norma{\cpt}_h+1}
\stackrel{j\to\infty}{\longrightarrow} 0
\end{eqnarray*}
since $\norma{\cpt}_h$ is a positive constant. By the arbitrariness of $y\in S_{\Pi_j}(\omega)$ we obtain the thesis.
\end{mynumerate}
\end{myproof}
\begin{myproof}[ of Proposition \ref{pro: independence of partition}]
Let $\Pi_j$ and $\Pi_l'$ be two partitions of the two distinct refinement
sequences $\graffa{\Pi_j}_{j\in\bbN}$ and $\graffa{\Pi_l'}_{l\in\bbN}$ of
$[t_0,t]$. Let $\Pi''=\Pi_j\cup\Pi_l'$ be the refinement of both $\Pi_j$
and $\Pi_l'$. Then Proposition \ref{pro: sn<sm and Sm<Sn} and Proposition
\ref{pro: sn<Sn} imply that $s_{\Pi_j}\subseteq s_{\Pi''} \subseteq
S_{\Pi''}\subseteq S_{\Pi_l'}$. Therefore $ s_{\Pi_j}\subseteq S_{\Pi_l'}$ for every $j,l\in\bbN$.
Then
$$
\overline{\tonda{\bigcup_{j\in\bbN} s_{\Pi_j}}} \subseteq
\bigcap_{l\in\bbN} S_{\Pi_l'}.
$$
Analogously
$$
\overline{\tonda{\bigcup_{l\in\bbN} s_{\Pi_l'}}} \subseteq
\bigcap_{j\in\bbN} S_{\Pi_j}.
$$
Proposition \ref{pro: dist(sn,Sn)--> 0} concludes the proof.
\end{myproof}

In order to prove Theorem \ref{teo: theta increasing}, let us consider the
following Lemma.
\begin{lem}\label{lem: s and S increasing}
Let $s,t\in [t_0,T]$ with $t_0<s<t$ and let $\Pi^s$ and $\Pi^t$ be two
partition of $[t_0,s]$ and $[t_0,t]$ respectively, such that $\Pi^s\subset
\Pi^t$. Then
$$
s_{\Pi^s}(s)\subseteq s_{\Pi^t}(t)\qquad \textrm{ and }\qquad
S_{\Pi^s}(s)\subseteq S_{\Pi^t}(t).
$$
\end{lem}
\begin{myproof}
The proofs of the two inclusions are similar. Let us prove that
$s_{\Pi^s}(s)\subseteq s_{\Pi^t}(t)$.
\\%
Since $\Pi^s\subset\Pi^t$, then $\Pi^s=\tonda{t_i}_{i=0}^n$ and
$\Pi^t=\Pi^s\cup\tonda{t_i}_{i=n+1}^{n+m}$ with $m\in\bbN$. By Lemma
\ref{lem: int_i G < int_I G}, we have that
\begin{eqnarray*}
s_{\Pi^s}(s) &=& \tonda{\nucleation_{t_0} \Mink \int_{t_0}^s \growth(\tau) d\tau
}\cup \bigcup_{i=1}^{n} \tonda{\Delta \nucleation_{t_i}\Mink \int_{t_i}^s
\growth(\tau) d\tau}
\\%
&\subseteq & \tonda{\nucleation_{t_0} \Mink \int_{t_0}^t \growth(\tau) d\tau
}\cup \bigcup_{i=1}^{n} \tonda{\Delta \nucleation_{t_i}\Mink \int_{t_i}^t
\growth(\tau) d\tau}
\\%
&\subseteq & \tonda{\nucleation_{t_0} \Mink \int_{t_0}^t \growth(\tau) d\tau
}\cup \bigcup_{i=1}^{n} \tonda{\Delta \nucleation_{t_i}\Mink \int_{t_i}^t
\growth(\tau) d\tau}
\\%
&& \cup \bigcup_{i=n+1}^{n+m} \tonda{\Delta \nucleation_{t_i}\Mink \int_{t_i}^t
\growth(\tau) d\tau}
\end{eqnarray*}
i.e.\ $s_{\Pi^s}(s)\subseteq s_{\Pi^t}(t)$.
\end{myproof}
\begin{myproof}[ of Theorem \ref{teo: theta increasing}]
For every $s,t\in [t_0,T]$ with $s<t$, let $\graffa{\Pi^s_i}_{i\in\bbN}$
and $\graffa{\Pi^t_i}_{i\in\bbN}$ be two refinement sequences of $[t_0,s]$
and $[t_0,t]$ respectively, such that $\Pi^s_i\subset \Pi^t_i$ for every
$i\in\bbN$. Then, by Lemma \ref{lem: s and S increasing}, $
S_{\Pi^s_i}\subseteq S_{\Pi^t_i} $. Now, as $i$ tends to infinity, we
obtain
$$
\Theta_s= \bigcap_{i\rightarrow \infty} S_{\Pi^s_i} \subseteq
\bigcap_{i\rightarrow \infty} S_{\Pi^t_i}= \Theta_t.
$$

For the second part, note that Theorem \ref{teo: any lebesg-integral of G
is RACS} still holds replacing $\salgebra_t$ instead
of $\salgebra$, so that for every $s\in[t_0,T]$, the family $\graffa{
\int_{s}^t \growth(\omega,\tau) d\tau}_{t\in [s,T]}$ is an adapted process
to the filtration $\graffa{\salgebra_t}_{t\in[t_0,T]}$. This fact together with Assumption \ref{hp(cont): Ht RACS} guarantees that $\graffa{S_{\Pi}}_{t\in[s,T]}$ is adapted for every partition $\Pi$ of $[s,T]$ and hence $\Theta$ is adapted too.
\end{myproof}

\section{Discrete time case and infinitesimal notations}%
\label{sec:discrete_time_case}

Let us consider $\Theta_s$ and $\Theta_{t}$ with $s< t$. Let $\graffa{\Pi^s_j}_{j\in\bbN}$ and $\graffa{\Pi^t_j}_{j\in\bbN}$ be two refinement sequences of $[t_0,s]$ and $[t_0,t]$ respectively, such that $\Pi^s_j\subset \Pi^t_j$ for every $j\in\bbN$ (i.e.\ $\Pi_j^s=\tonda{t_i}_{i=0}^n$ and $\Pi_j^t=\Pi_j^s\cup\tonda{t_i}_{i=n+1}^{n+m}$ with $n,m\in\bbN$).
It is easy to compute
$$
s_{\Pi^t_j} = \tonda{s_{\Pi^s_j} \Mink \int_{s}^t \growth(\tau) d\tau }
\cup \bigcup_{i=n+1}^{n+m} \tonda{\Delta \nucleation_{t_i}\Mink \int_{t_i}^t
\growth(\tau) d\tau}.
$$
Then, by Definition \ref{def: continuous set process}, whenever $\modulo{\Pi^t_j}\to 0$, we obtain
\begin{equation}\label{eq: continuous like discrete}
\Theta_t = \tonda{\Theta_s \Mink \int_{s}^t \growth(\tau) d\tau } \cup
\lim_{\modulo{\Pi^t_j}\rightarrow 0}\bigcup_{i=n+1}^{n+m} \tonda{\Delta
\nucleation_{t_i}\Mink \int_{t_i}^t \growth(\tau) d\tau}.
\end{equation}
The following notations
$$
\growth_k= \int_{s}^t \growth(\tau) d\tau \nonumber
\qquad\textrm{and}\qquad
\nucleation_k = \lim_{\modulo{\Pi^t_j}\rightarrow 0}\bigcup_{i=n+1}^{n+m} \tonda{\Delta \nucleation_{t_i}\Mink \int_{t_i}^t \growth(\tau) d\tau}\nonumber
$$
lead us to the set-valued discrete time stochastic process
$$
\Theta_k=\left\{
\begin{array}{ll}
(\Theta_{k-1}\Mink \growth_{k} ) \cup \nucleation_k, & k\ge 1,\\
\nucleation_0, & k=0.
\end{array}\right.
$$
In view of this, we are able to justify infinitesimal notations introduced in \eqref{eq: infinitesimal continuous like discrete}. In particular, from Equation \eqref{eq: continuous like discrete}, whenever $\modulo{\Pi^t_j}\rightarrow 0$, we obtain
$$
\Theta_t=
\tonda{\nucleation_{t_0} \Mink \int_{t_0}^t \growth(\tau) d\tau }\cup
\bigcup_{s=t_0}^t \tonda{d\nucleation_s \Mink \int_{s}^t \growth(\tau) d\tau}, 
\qquad t\in [t_0,T].
$$
Moreover, with a little abuse of this infinitesimal notation, we get two differential formulations
$$
d\Theta_t = \oplus \growth_{t} dt \cup d\nucleation_{t}
\qquad\textrm{and}\qquad \Theta_{t+dt}=(\Theta_t\Mink \growth_{t}
dt)\cup d\nucleation_{t}.
$$

%
%
%
%

\end{document}